%% file: gen_goedel_all1.tex
\documentclass[11pt]{article}
\usepackage{latexsym,a4}

\usepackage[cp1250]{inputenc}
\usepackage[russian,english]{babel}
\usepackage{amsmath,amsfonts,amsthm,amssymb}
\usepackage{bussproofs}

\EnableBpAbbreviations

\input{macros.tex}

\begin{document}

\title{Some abstract versions of G\"odel's second incompleteness theorem based on non-classical logics}

\author{Lev Beklemishev and Daniyar Shamkanov \\
Steklov Mathematical Institute, Moscow\\
Moscow M.V.~Lomonosov State University \\
National Research University Higher School of Economics, Moscow\\
}

\maketitle

{\em To Albert Visser, a remarkable logician and a dear friend, whose papers and conversations are a source of constant inspiration}

\begin{abstract} We study abstract versions of G\"odel's second incompleteness theorem and formulate generalizations of L\"ob's derivability conditions that work for logics weaker than the classical one. We isolate the role of  contraction rule in G\"odel's theorem and give a (toy) example of a system based on modal logic without contraction invalidating G\"odel's argument.
\end{abstract}

\section{Introduction}
One of the topics that have been fascinating logicians over the years is G\"odel's second incompleteness theorem (G2).
Both mathematically and philosophically G2 is well known to be more problematic than his first incompleteness theorem (G1).
G1 and Rosser's Theorem are well understood in the context of recursion theory. Abstract logic-free formulations have been given by Kleene \cite{Kle50} (`symmetric form'), Smullyan \cite{Smu94} (`representation systems') and others. Sometimes G2 is considered as a minor addition to G1, whose role is to exhibit a specific form of the sentence independent from a given theory, namely its consistency assertion. However, starting with the work of Kreisel, Orey, Feferman, and others, who provided various nontrivial uses of G2, it has been gradually understood that the two results are of a rather different nature and scope. G2 has more to do with the (modal-logical) properties of the provability predicate and the phenomenon of self-reference in sufficiently expressive systems. A satisfactory general mathematical context for G2, however, still seems to be lacking.

The main difficulties in G2 are due to the fact that we cannot easily delineate a class of formulas that `mean' consistency. Thus, the most intuitively appealing formulation of G2 --- \emph{sufficiently strong consistent theories cannot prove their own consistency} --- remains non-mathematical.
For a concrete formal system, such as Peano arithmetic $\PA$, one can usually write out a specific `natural' formula $\Con_\PA$ and declare it to be the expression of consistency. This approach is rather common in mathematics but has several deficiencies: Firstly, it ties the statement to a very particular formula, coding mechanism etc., and provides no clue why this choice is better than the other ones. Secondly, instead of a general theorem working uniformly for a wide class of theories, we only obtain a specific statement for an individual theory such as $\PA$. We do not know what is the natural consistency assertion for an arbitrary extension of $\PA$. Thus, we have a problem with translating our informal intuition into strict mathematical terms.

 The way to better understand G2 is through investigating its range and generalizations. A lucky circumstance is that G2 also holds for larger syntactically defined classes of consistency formulas, some of which are apparently intensionally correct (adequately express consistency), but some are not. Thus, it is still possible to formulate mathematical results in certain important aspects
\emph{more} (rather than less) general than the broad intuitive formulation of G2 above.

A universally accepted approach to general formulations of G2 appeared in the fundamental paper by Feferman~\cite{Fef60} who showed, among other things, that G2 holds for all consistency assertions defined by $\Sigma_1$-numerations.
Feferman deals with first-order theories $T$ in the language containing that of $\PA$ and specified by recursively enumerable (r.e.)\ sets of axioms. Feferman assumes fixed some natural G\"odel numbering of the syntax of $T$ as well as some specific axiomatization of first order logic. A $\Sigma_1$-formula $\alpha(x)$ defining the set of G\"odel numbers of axioms of $T$ in the standard model of $\PA$ is called a \emph{$\Sigma_1$-numeration} of $T$.\footnote{Feferman deals with the notion of r.e.\ formula rather than with the equivalent notion of $\Sigma_1$-formula more common today.} It determines the
provability formula $\Prs{\ga}(x)$ and the corresponding consistency assertion $\Con_\ga$. Feferman's statement of G2 is that for all consistent theories $T$ given by $\Sigma_1$-numerations $\ga$ and containing a sufficiently strong fragment of $\PA$, the formula $\Con_\ga$ is unprovable in $T$.

This theorem is considerably more general than any specific instance of G2 for an individual theory $T$. However, it also presupposes quite a lot: first order logic and its axiomatization, G\"odel numbering, the way formula $\Prs{\ga}$ is built from $\ga$.

Exploring bounds to G2 leads to relaxing various assumptions involved in Feferman's statement:

\bi
\item One can weaken the axioms of arithmetic (for a representative selection see Bezboruah--Shepherdson~\cite{BezShe}, Pudl\'ak~\cite{Pud85}, Wilkie--Paris~\cite{WP}, Adamowicz--Zdanowski~\cite{AZ11}, Willard~\cite{Wil01,Wil06}).
\item One can consider theories modulo interpretability. This approach started with the work of Feferman~\cite{Fef60}. In the recent years it lead to particularly attractive coding-free formulations of generalizations of G2 due to Harvey Friedman and Albert Visser~(see \cite{Vis93,Vis11,Vis12}).
\item One can weaken the requirements on the proof predicate aka derivability conditions (see Feferman~\cite{Fef60}, L\"ob~\cite{Lob55}, Jeroslow~\cite{Jer70,Jer73}).
\item One can weaken the logic.
\ei

It is the latter two aspects, less studied in the literature, that we are going to comment on in this note. Firstly, let us briefly recall the history of derivability conditions.

G\"odel \cite{God31} gave a sketch of a proof of G2 and a promise to provide full details in a subsequent publication. This promise has not been fulfilled, and a detailed proof of this theorem --- for a system $Z$ related to first-order arithmetic $\PA$ --- only appeared in a monograph by Hilbert and Bernays \cite{HB}. In order to structure a rather lengthy proof Hilbert and Bernays formulated certain conditions on the proof predicate in $Z$, sufficient for the proof of G2. Later Martin L\"ob~\cite{Lob55} gave an elegant form to these conditions by stating them fully in terms of the provability predicate $\Pr(x)$ and obtained an important strengthening of G2 known as L\"ob's Theorem. Essentially the same properties of the provability predicate were earlier noted by G\"odel in his note~\cite{God33}, where he proposed to treat the provability predicate as a connective $\Box$ in modal logic, though the idea that these conditions constitute necessary requirements on a provability predicate most likely only appeared later. For the sake of brevity we call the G\"odel--Hilbert--Bernays--L\"ob conditions simply \emph{L\"ob's conditions} below.

A traditional proof of G2 (for arithmetical theories) consists of a derivation of G2 from the fixed point lemma using L\"ob's conditions (see e.g.~\cite{Smo77}). An accurate justification of these conditions is technically not so easy, and a rare textbook provides enough details here, however see Smory\'nski~\cite{Smo85} and Rautenberg~\cite{Rau06} for readable expositions.

L\"ob's conditions are applicable to formal theories at least containing the connective of implication and closed under the \emph{modus ponens} rule. Here we give more general abstract formulations of G2 which presuppose very little about logic. They are rather close in the spirit and the level of generality to the recursion-theoretic formulations of G1 due to Smullyan. When a good implication is added to the language one essentially obtains the familiar L\"ob's conditions. However, we show that G\"odel's argument presupposes admissibility of the contraction rule restricted to $\Box$-formulas in the logic under consideration. Moreover, the uniqueness of G\"odelian fixed point is based on the similarly restricted form of weakening.

In the last part of the paper we present a system invalidating a formalized version of G2. We consider a version of propositional modal logic K4 based on the contraction-free fragment of classical logic extended by fixed point operators (defined for any formulas modalized in the fixed point variables). By means of a cut-elimination theorem for this system we establish the failure of G2 and some other properties such as the infinity of the  G\"odelian and Henkinian fixed points.

\section{Abstract provability structures}

\bd \label{cns}
Let us call \emph{an abstract consequence relation} a structure $S=(L_S,\leq_S,\top,\bot)$, where $L_S$ is a set of  \emph{sentences} of $S$, $\leq_S$ is a transitive reflexive relation on $L_S$, $\top$ and  $\bot$ are distinguished elements of $L_S$ (`axiom' and `contradiction'). A sentence $x\in L_S$ is called \emph{provable in $S$}, if $\top\leq_S x$, and \emph{refutable in $S$}, if $x\leq_S \bot$. Sentences $x,y$ are called \emph{equivalent in $S$}, if $x\leq_S y$ and $y\leq_S x$. The equivalence of $x$ and $y$ will be denoted $x=_S y$.
\ed

The structure $S$ represents syntactical (rather than semantical) data about the theory in question.
In a typical case, for example, for arithmetical theories $S$, the relation $x\leq_S y$ denotes the provability of $y$ from hypothesis $x$, whereas $\top$ and $\bot$ are some standard provable and refutable formulas, respectively, e.g., $0=0$ and $0\neq 0$.

In concrete situations we can enrich this structure by additional data, for example, by the conjunction and the implication connectives. Notice that we do not assume either $\bot\leq_S x$ or $x\leq_S \top$, nor do we assume the existence of any logical connectives (such as negation) in $S$.

$S$ is called \emph{inconsistent} if $\top\leq_S \bot$, otherwise it is called \emph{consistent}. By transitivity, if $S$ is consistent then no sentence is both provable and refutable. $S$ is called \emph{complete} if every $x\in L_S$ is either provable or refutable. $S$ is called \emph{r.e.}, if $L_S$ is recursive and $\leq_S$ is r.e.\ (as a binary relation). $T$ is called an \emph{extension} of $S$ if $L_T=L_S$ and $\leq_S$ is contained in $\leq_T$.

Let $P_S$ and $R_S$ denote the sets of provable and of refutable sentences of $S$, respectively. If $S$ is consistent and r.e., then $P_S$ and $R_S$ is a pair of disjoint r.e.\ sets. We say that $S$ \emph{separates pairs of disjoint r.e.\ sets} if for each such pair $(A,B)$ there is a total computable function $f$ such that
$$\al{n\in A} f(n)\in P_S \text{ and } \al{n\in B} f(n)\in R_S.$$

The following statement is a natural version of G1 and Rosser's theorem for abstract consequence relations (\'a la Kleene and Smullyan); we omit the standard proof.
\bpr
\begin{enumr}
\item If $S$ is r.e., consistent and complete, then both $P_S$ and $R_S$ are decidable.
\item If $S$ is r.e.\ and separates disjoint pairs of r.e.\ sets, then every consistent extension of $S$ is incomplete and undecidable.
\end{enumr}
\epr

Next we introduce two operators $\Box,\boxtimes:L_S\to L_S$ representing provability and refutability predicates in $S$.
\bd
\emph{Provability} and \emph{refutability operators} for an abstract consequence relation $S$ are functions $\Box,\boxtimes:L_S\to L_S$ satisfying the following conditions, for all $x,y\in L_S$:
\label{prop-def}
\renewcommand{\labelenumi}{C\arabic{enumi}.}
\ben
\item $x\leq_S y \ \Imp \ \Box x\leq_S \Box y$, $\boxtimes y\leq \Boxt x$. \label{g1}
\item $\top\leq_S \Boxt \bot$; \label{g2}
\item $x\leq_S\Box y,\ x\leq_S\Boxt y\ \Imp x\leq_S \Boxt \top$; \label{g3}
\item $\Boxt x\leq_S\Box\Boxt x$. \label{g4} 
\een

The algebra $(L_S,\leq_S,\top,\bot,\square,\Boxt)$ is called an \emph{abstract provability structure} (APS).
\ed

Intuitively, $\Box x$ is the sentence expressing the provability of a sentence $x$, whereas $\Boxt x$ expresses its refutability in $S$. Condition C\ref{g1} means that provability of $y$ follows from provability of $x$ whenever $y$ is derivable from $x$; similarly, refutability of $y$ implies refutability of $x$. Conditions C\ref{g2} and C\ref{g3} are axioms for contradiction: according to C\ref{g2}, refutability of $\bot$ is provable in $S$; according to
C\ref{g3}, $\top$ is refutable if some sentence $y$ is both provable and refutable.
Finally, Condition \ref{g4} means that the refutability of $x$ can be formally checked in $S$. It is an analogue of L\"ob's condition L2 (see below).

Note that we consider the refutability operator on a par with the provability operator, since we do not assume that the logic of $S$ necessarily has a well-defined operation of negation, that is, we cannot always define $\Boxt x$ as $\Box\neg x$.

\brem It is rather natural to additionally require that $\Box\bot=_S\Boxt\top$: refutability of $\top$ and provability of $\bot$ are expressed by the same statement $\top\leq_S\bot$. Yet, it is not, strictly speaking, needed in this very abstract context, and we take $\Boxt\top$ as our default expression of inconsistency.
\erem

\bd We say that an abstract provability structure $S$ \emph{has a G\"odelian fixed point} if there is a sentence $p\in L_S$ such that $p=_S\Boxt p$.
\ed

Notice that G\"odel considered a dual sentence $q$ expressing its own unprovability in $S$. R. Jeroslow \cite{Jer73} noticed that the sentence stating its own refutability allows to prove G2 under somewhat more general conditions than those of L\"ob. In our formalism the sentence $q$ is not expressible, therefore we are using Jeroslow's idea.

A very abstract version of G2 can now be stated as follows.

\bt \label{god2-abs}
Suppose an APS $S$ has a G\"odelian fixed point.
\begin{enumr}
\item If $S$ is consistent, then $\Boxt \top$ is irrefutable in $S$.
\item $\Boxt\Boxt \top \leq_S \Boxt \top$, that is, Statement \emph{(i)} is formalizable in $S$.
\end{enumr}
\et

\bp\ Let $p=_S\Boxt p$. First we prove Statement (ii) omitting the subscript $_S$ everywhere:

\ben
\item $\Boxt p\leq \Box\Boxt p\leq \Box p$ by C\ref{g4} and C\ref{g1};
\item $p=\Boxt p\leq \Boxt\top$ by C\ref{g3} (since $\Boxt p\leq \Boxt p$);
\item $\Boxt\Boxt\top\leq \Boxt p=p\leq\Boxt \top$ by C\ref{g1}.
\een

Proof of Statement (i):
Assume $\Boxt \top\leq \bot$. By the previous argument $p\leq \Boxt \top$, hence $p\leq\bot$. By C\ref{g1}, $\Boxt \bot\leq \Boxt p = p \leq \bot$.
Therefore, by C\ref{g2}, $\top\leq \Boxt \bot\leq \bot$. \ep

The following statement shows that under some additional condition the G\"odelian--Jeroslowian fixed point is unique modulo equivalence in $S$ and coincides with the inconsistency assertion for $S$. Therefore, the existence of such a fixed point is not only sufficient but also necessary for the validity of (a formalized version of) G2. The additional condition is
\bi
\item[C5.] $x\leq_S \top$, for all $x\in L_S$.
\ei

\bt \label{fixp} Assume C5 holds for $S$. Then $p=_S\Boxt \top$ for all G\"odelian fixed points $p$ and (if such a sentence exists) $$\Boxt\Boxt \top =_S \Boxt \top.$$
\et

\bp\ We know that $p\leq \Boxt \top$. Since $p=\Boxt p\leq\top$ we obtain  $\Boxt\top\leq \Boxt p = p$. Hence $p=\Boxt\top$ and therefore  $\Boxt\Boxt \top = \Boxt \top.$
\ep

\section{Consequence relations with implication}

Classical L\"ob's conditions emerge for APS with an implication. A decent implication can be defined for consequence relations representing derivability of a sentence from a (multi)set of assumptions. In other words, we now go to a more general but less symmetric format $\Gamma\vdash\phi$, where $\Gamma$ is a finite multiset and $\phi$ an element of a given set $L$. In order to avoid confusion we use the more standard notation $\vdash$ instead of $\leq$ and will follow the standard conventions of sequential proof format. In particular, $\Gamma,\phi$ denotes the result of adjoining $\phi\in L_S$ to a multiset of sentences $\Gamma$, and $\Gamma,\Delta$ denotes the multiset union of $\Gamma$ and $\Delta$.\footnote{Our strive for generality does not go as far as to consider lists of formulas rather than multisets.}

\bd \emph{A consequence relation with an implication on $L$} is a structure $S=(L_S,\vdash,\to,\top,\bot)$ where
$\vdash$ is a binary relation between finite multisets of elements of $L_S$ and elements of $L_S$; $\to$ is a binary operation on $L$; $\top$ and $\bot$ are distinguished elements of $L$ such that the following conditions hold:

\ben
\renewcommand{\labelenumi}{I\arabic{enumi}.}
\item $\phi\vdash \phi$;
\item if $\Gamma,\psi\vdash \phi$ and $\Delta\vdash\psi$ then $\Gamma,\Delta\vdash\phi$;
\item $\Gamma,\phi\vdash\psi \iff \Gamma\vdash \phi\to\psi$;
\item $\Gamma,\top\vdash\phi\iff \Gamma\vdash\phi.$
\een
\ed

Notice that Conditions I1 and I2 generalize reflexivity and transitivity of $\leq$. Setting $\phi\leq_S\psi$ as $\phi\vdash\psi$ yields an abstract consequence relation in the sense of Definition~\ref{cns}. Condition I3 speaks for itself. Condition I4 conveniently stipulates that provability from the empty multiset of assumptions is the same as provability from $\top$. It also implies $\top\to\bot=_S\bot$.

Similarly to the implication one can consider consequence relations with other additional connectives of which we are mostly interested in conjunction.

\bd \emph{Conjunction} is a binary operator $\otimes:L_S^2\to L_S$ satisfying
$$\Gamma,\phi,\psi\vdash \theta \iff \Gamma, \phi\otimes\psi\vdash\theta.$$
\ed

If conjunction is available, then $\phi_1,\dots,\phi_n\vdash\psi$ holds in $S$ if and only if  $\phi_1\otimes\cdots\otimes\phi_n\vdash\psi$. Hence, in the presence of conjunction in $S$ the relation $\leq_S$ uniquely determines the corresponding multiset consequence relation $\Gamma\vdash\phi$.

For consequence relations with an implication we can define negation $\neg\phi$ by $\phi\to\bot$. The following simple lemma shows that the implication respects the deductive equivalence relation in $S$ and the negation satisfies the contraposition principle.

\bl \label{imp-prop}
\begin{enumr}
\item If $\Gamma\vdash \phi\to\psi$ and $\Delta\vdash\phi$, then $\Gamma,\Delta\vdash\psi$;
\item  $\phi_1=_S\phi_2$ and $\psi_1=_S\psi_2$ implies $ (\phi_1\to\psi_1)=_S (\phi_2\to\psi_2);$
\item $\Gamma,\phi\vdash\psi$ implies $\Gamma,\neg\psi\vdash\neg\phi$.
\end{enumr}
\el

Next we turn to the derivability conditions. Assume $S$ is a consequence relation with an implication.

\bd $\Box:L_S\to L_S$ satisfies L\"ob's derivability conditions for $S$ if
\ben
\renewcommand{\labelenumi}{L\arabic{enumi}.}
\item $\Box(\phi\to\psi)\vdash \Box\phi\to\Box\psi$;
\item $\Box\phi\vdash\Box\Box\phi$;
\item $\vdash\phi$ implies $\vdash\Box\phi$.
\een
\ed

\bl For any consequence relation with an implication the following statements are equivalent:
\begin{enumr}
\item $\Box$ satisfies L\"ob's conditions for $S$;
\item $\Box$ satisfies L2 and $S$ is closed under the rule
$$\frac{\Gamma\vdash\phi}{\Box\Gamma\vdash\Box\phi};$$
\item $S$ is closed under the rule
$$\frac{\Gamma,\Box\Delta\vdash\phi}{\Box\Gamma,\Box\Delta\vdash\Box\phi}.$$
\end{enumr}
\el

\brem Notice that the last rule is formulated slightly differently from the more standard rule for modal logic K4:
$$\frac{\Gamma,\Box\Gamma\vdash\phi}{\Box\Gamma\vdash\Box\phi}.$$
The latter has a form of built-in contraction that we are not assuming here.
\erem

It is natural to define refutability $\boxtimes\phi$ as provability of negation $\Box\neg\phi$. Notice that since $\bot=_S\top\to\bot$ we have $\Boxt\top=_S \Box\bot$, whenever L1 holds for $\Box$. However, as the example in Section \ref{sec-ex1} shows, this translation does not always yield an APS in the sense of Definition \ref{prop-def}. To sort things out we need to consider two additional conditions on the consequence relation.

\bd A consequence relation with an implication
\bi
\item[-] \emph{satisfies contraction} if
$\Gamma,\phi,\phi\vdash\psi$ implies $\Gamma,\phi\vdash\psi$;
\item[-] \emph{satisfies weakening} if
$\Gamma\vdash\psi$ implies $\Gamma,\phi\vdash\psi$, for any $\phi$.
\ei
\ed

The first condition intuitively means that any hypothesis can be used several times in a derivation. Recall that for Girard's linear logic this condition is not met, however it is postulated, for example, for relevant logics. It turns out that a certain amount of contraction is essential for the proof of G2.

The second condition corresponds to the requirement $x\leq_S \top$ that was needed to guarantee that $\Boxt \top$ is a G\"odelian fixed point and that such a fixed point is unique.

For consequence relations with an implication we have the following proposition.

\bpr \label{s-lob}
Suppose $S$ satisfies contraction, $\Box:L_S\to L_S$ satisfies L\"ob's conditions for $S$ and $\Boxt\phi:=\Box(\phi\to\bot)$. Then $(L_S,\leq_S,\Box,\boxtimes,\top,\bot)$ is an APS.
\epr

\bp By Lemma \ref{imp-prop} $\phi\vdash\psi$ implies $\neg\psi\vdash\neg\phi$. This yields Conditions C\ref{g1} and C\ref{g4}. Condition C\ref{g2} obviously follows from Condition 1 for a good consequence relation. Let us prove C\ref{g3}.
By Lemma \ref{imp-prop}(i) we have: $\phi,\neg\phi,\top\vdash \bot$. Hence, $\phi,\neg\phi\vdash\top\to\bot$, therefore $\Box\phi,\Box\neg\phi\vdash\Box\neg\top$ by Condition 1. The rules of transitivity and contraction imply that, if  $\Gamma\vdash\Box\phi$ and $\Gamma\vdash\Box\neg\phi$, then $\Gamma\vdash\Box\neg\top$. \ep

Thus, from Proposition \ref{s-lob} we obtain the following expected corollary, parallel to Theorem \ref{god2-abs}, for consequence relations satisfying contraction.

\bt \label{g2contr} Suppose $S$ satisfies contraction and $\Box$ satisfies L\"ob's conditions for $S$.
Then Theorem \ref{god2-abs} holds for $S$.
\et

For an analogue of Theorem \ref{fixp} on the uniqueness of a G\"odelian fixed point we also need a weakening property.

\bt \label{g2weak} Suppose $S$ satisfies contraction and weakening and $\Box$ satisfies L\"ob's conditions for $S$. Then all G\"odelian fixed points in $S$ (if exist) are equivalent to $\Boxt\top=_S\Box\bot$. \et

\brem
As it turns out, contraction and weakening for $S$, though natural, are somewhat excessive requirements for the validity of Theorems \ref{g2contr} and \ref{g2weak}. A consequence relation with an implication
\bi
\item[-] \emph{satisfies $\Box$-contraction} if
$\Gamma,\Box\phi,\Box\phi\vdash\psi$ implies $\Gamma,\Box\phi\vdash\psi$;
\item[-] \emph{satisfies $\Box$-weakening} if
$\Gamma\vdash\phi$ implies $\Gamma,\Box\psi\vdash\phi$, for any $\psi$.
\ei

Conditions C3 and C5 of APS can also be weakened to
\bi
\item[C3$'$.] $\Boxt x\leq_S\Box y,\ \Boxt x\leq_S\Boxt y\ \Imp \Boxt x\leq_S \Boxt \top$;
\item[C5$'$.] $\Boxt x\leq_S \top$.
\ei
With these modifications, the proofs of Theorems \ref{god2-abs} and \ref{fixp} stay the same, which in turn yields more general versions of Theorems \ref{g2contr} and \ref{g2weak} for consequence relations satisfying only $\Box$-contraction and $\Box$-weakening.

The property of $\Box$-contraction actually holds for some meaningful arithmetical systems lacking general contraction rule, for example, for a version of Peano arithmetic based on affine predicate logic considered by the second author of this paper (as yet, unpublished).
\erem

\section{A non-G\"{o}delian theory with fixed points} \label{sec-ex1}
In view of Theorems \ref{g2contr} and \ref{g2weak} it is natural to ask whether the assumptions of $\Box$-contraction and $\Box$-weakening are substantial for these results. More specifically, two questions immediately present themselves:
\ben
\item Does there exist a consequence relation with an implication satisfying L\"ob's conditions for $\Box$ in which a G\"odelian fixed point exists, but G2 fails? (The failure of G2 can be understood in two different senses --- as a failure of its formalized version, and as a failure of its non-formalized version. Our example will show the failure of the formalized version.)
\item Do G\"odelian fixed points in such a system $S$ have to be unique, even if $S$ satisfies weakening?
\een

In this section we provide an example showing that the answer to the first question is positive and to the second one negative. Moreover, we formulate a system in which there are many more fixed points than are officially required for a proof of G2. Our system $\mathsf{S}$ is a version of modal logic K4 based on the multiplicative $\{\to,\otimes,\bot\}$ fragment of a classical logic without contraction. It also has a built-in fixed point operator where the expression $\mathsf{fp}\, x. A(x)$ denotes some fixed point of $A(x)$ for formulas $A$ modalized in the variable $x$. Thus, one will be able to derive $$\mathsf{fp}\, x. A(x)=_S A(\mathsf{fp}\, x. A(x)),$$
for each formula $A(x)$ modalized in $x$.
Let us now turn to the exact definitions.

\medskip
Consider the set of formulas $\mathsf{Fm}_0$ given by the grammar:
$$ A ::= p \,\,|\,\, x \,\,|\,\, \bot \,\,|\,\, (A \rightarrow A) \,\,|\,\,  \Box A \;, $$
where $p$ stands for \emph{atomic propositions} and $x$ stands for \emph{variables} (the alphabets of atomic propositions and variables are disjoint). We define the \emph{set of formulas of $\mathsf{S}$} by extending the set $\mathsf{Fm}_0$ by a new constructor: if $A$ is a formula and all free occurrences of $x$ in $A$ are within the scope of modal operators, then $\mathsf{fp}\, x. A$ is a formula, and $\mathsf{fp}\, x$ binds all free occurrences of $x$. A formula $B$ is \emph{closed} if it does not contain any free occurrences of variables. For a closed formula $B$, we denote by $A[B/\!/x]$ the result of replacing all free occurrences of $x$ in $A$ by $B$. 
We also put $\neg A  := A\rightarrow \bot$, $\top  := \neg \bot$ and $A \otimes B  := \neg (A \rightarrow \neg B)$.

A \textit{sequent} is an expression of the form $\Gamma \Rightarrow \Delta$, where $\Gamma$ and $\Delta$ are finite multisets of closed formulas. 
The sequent calculus $\mathsf{S}$ is defined in the standard way by the following initial sequents and inference rules:

\begin{gather*}
\AXC{ $\Gamma, A \Rightarrow A, \Delta $}
\DisplayProof \qquad
\AXC{ $\Gamma , \bot \Rightarrow  \Delta$}
\DisplayProof \\\\
\AXC{$\Gamma, A[\mathsf{fp}\, x.\, A/\!/x] \Rightarrow \Delta$}
\LeftLabel{$(\mathsf{fix_L})$}
\UIC{$\Gamma, \mathsf{fp}\, x.\, A \Rightarrow  \Delta$}
\DisplayProof \qquad
\AXC{$\Gamma \Rightarrow  A[\mathsf{fp}\, x.\, A/\!/x],\Delta$}
\LeftLabel{$(\mathsf{fix_R})$}
\UIC{$\Gamma \Rightarrow  \mathsf{fp}\, x.\, A ,\Delta$}
\DisplayProof
\end{gather*}
\begin{gather*}
\AXC{$\Gamma , B \Rightarrow  \Delta$}
\AXC{$\Sigma \Rightarrow  A, \Pi$}
\LeftLabel{$(\mathsf{\rightarrow_L})$}
\BIC{$\Gamma , \Sigma, A \rightarrow B \Rightarrow  \Pi, \Delta$}
\DisplayProof \qquad
\AXC{$\Gamma, A \Rightarrow  B ,\Delta$}
\LeftLabel{$(\mathsf{\rightarrow_R})$}
\UIC{$\Gamma \Rightarrow  A \rightarrow B ,\Delta$}
\DisplayProof \\\\
\AXC{$\Sigma, \Box \Pi \Rightarrow A$}
\LeftLabel{$(\mathsf{\Box})$}
\UIC{$\Gamma, \Box \Sigma, \Box \Pi \Rightarrow \Box A , \Delta$}
\DisplayProof \;.
\end{gather*}
Explicitly displayed formulas in the conclusions of the rules are called \emph{principal formulas} of the corresponding inferences.
In the rules $(\mathsf{fix_L})$, $(\mathsf{fix_R})$, $(\mathsf{\rightarrow_L})$ and $(\mathsf{\rightarrow_R})$, the elements of
$\Gamma$, $\Delta$, $\Sigma$ and $\Pi$ are called \emph{side formulas}.
In initial sequents and in applications of the rule $(\mathsf{\Box})$, the elements of
$\Gamma$ and $\Delta$ are \emph{weakening formulas}. We call the elements of $\Box \Sigma$ and $\Box \Pi$ in the corresponding applications of $(\mathsf{\Box})$ \emph{active formulas}. In addition, explicitly displayed formulas in initial sequents are called \emph{axiomatic formulas}.

A \emph{proof in $\mathsf{S}$} is a finite tree whose nodes are marked by sequents and leaves are marked by initial sequents that is constructed according to the rules of the sequent calculus. A sequent $\Gamma \Rightarrow \Delta$ \emph{is provable in $\mathsf{S}$} if
there is a proof with the root marked by $\Gamma \Rightarrow \Delta$.

We associate with $\mathsf{S}$ a consequence relation with an implication and conjunction in the usual way by letting $\Gamma\vdash_\mathsf{S}\phi$ iff $\Gamma\Imp\phi$ is provable in $\mathsf{S}$. The main thing we need to prove about $\mathsf{S}$ is the closure of
$\mathsf{S}$ under the cut rule, which would show that $\Gamma\vdash_\mathsf{S}\phi$ is indeed a well-defined consequence relation (see Theorem~\ref{cutelim} below).

Since $\mathsf{S}$ is cut-free, the following propositions are easy to establish. Firstly, we obtain the failure of formalized G2.

\bpr
The sequent $\Box (\Box \bot \rightarrow \bot) \Rightarrow \Box \bot$ is not provable in $\mathsf{S}$.
\epr

Recall that an inference rule is called admissible (for a given proof system) if, for every instance of the rule, the conclusion is provable whenever all premises are provable.

\bpr
The L\"{o}b rule and the Henkin rule
\begin{gather*}
\AXC{$\Box A \Rightarrow A$}
\LeftLabel{$(\mathsf{L\ddot{o}b})$}
\RightLabel{ }
\UIC{$ \quad  \Rightarrow A$}
\DisplayProof
\qquad
\AXC{$\Box A \Rightarrow A$}
\AXC{$A\Rightarrow \Box A$}
\LeftLabel{$(\mathsf{Hen})$}
\RightLabel{ }
\BIC{$ \quad  \Rightarrow A$}
\DisplayProof
\end{gather*}
are not admissible in $\mathsf{S}$.
\epr
\begin{proof}
Consider the Henkin fixed point $\mathsf{fp}\, x.\, \Box x$. The sequent $\Rightarrow \mathsf{fp}\, x.\, \Box x$ is not provable in $\mathsf{S}$. Hence, the Henkin rule is not admissible and so is the stronger L\"ob rule.
\end{proof}

\bpr
There are infinitely many Henkinian and G\"odelian fixed points in $\mathsf S$.
\epr
\bp The routine of bound variables in $\mathsf S$ is such that the formulas
$\mathsf{fp}\, x_i.\Box x_i$ for graphically distinct variables $x_i$ are all inequivalent. (There is no rule of bound variables renaming and, in fact, it is easy to convince oneself that there are no cut-free proofs in $\mathsf S$ of the sequents $\mathsf{fp}\, x_i.\Box x_i\Imp \mathsf{fp}\, x_j.\Box x_j$, for $i\neq j$.) The same holds for the G\"odelian fixed points of $\mathsf S$.
\ep

\section{Cut-admissibility for $\mathsf S$}
For a proof of the cut-admissibility theorem for $\mathsf S$ we need the following standard lemma. Let the \emph{size $\lVert \pi \rVert$ of a proof $\pi$} be the number of nodes in $\pi$.

\bl
The weakening rule
\begin{gather*}
\AXC{$\Gamma \Rightarrow \Delta$}
\LeftLabel{$(\mathsf{weak})$}
\RightLabel{ }
\UIC{$ \Sigma , \Gamma \Rightarrow   \Delta, \Pi$}
\DisplayProof
\end{gather*}
is admissible for $\mathsf{S}$, and its conclusion has a proof of at most the same size as the premise.
\el

\bt \label{cutelim}
The cut rule
\begin{gather*}
\AXC{$\Gamma \Rightarrow \Delta, A$}
\AXC{$A, \Sigma \Rightarrow \Pi$}
\LeftLabel{$(\mathsf{cut})$}
\RightLabel{ ,}
\BIC{$\Gamma, \Sigma \Rightarrow  \Pi, \Delta$}
\DisplayProof
\end{gather*}
is admissible for $\mathsf{S}$. Moreover, if $\pi_1$ and $\pi_2$ are proofs of the premises of $(\mathsf{cut})$, then the conclusion of $(\mathsf{cut})$ has a proof with the size being less than $ \lVert \pi_1 \rVert + \lVert \pi_2 \rVert$.
\et
\begin{proof}

Assume we have an inference
\begin{gather*}
\AXC{$\pi_1$}
\noLine
\UIC{\vdots}
\noLine
\UIC{$\Gamma \Rightarrow \Delta, A$}
\AXC{$\pi_2$}
\noLine
\UIC{\vdots}
\noLine
\UIC{$A, \Sigma \Rightarrow \Pi$}
\LeftLabel{$(\mathsf{cut})$}
\RightLabel{ ,}
\BIC{$\Gamma, \Sigma \Rightarrow  \Pi, \Delta$}
\DisplayProof
\end{gather*}
where $\pi_1$ and $\pi_2$ are proofs in $\mathsf{S}$. We proof by induction on $\lVert \pi_1 \rVert + \lVert \pi_2 \rVert$ that for any formula $A$ there exists a proof $\mathcal{E}_A(\pi_1,\pi_2)$ of $\Gamma, \Sigma \Rightarrow  \Pi, \Delta$ with the size being less than $\lVert \pi_1 \rVert + \lVert \pi_2 \rVert $.

Consider the final inference in $\pi_1$. If the formula $A$ is in a position of a weakening formula in it, then we erase $A$ in $\pi_1$ and extend the sequent $\Gamma \Rightarrow \Delta$ to $\Gamma, \Sigma \Rightarrow  \Pi, \Delta$ by adding new weakening formulas. This transformation of $\pi_1$ defines $\mathcal{E}_A(\pi_1,\pi_2)$. Moreover, we have $\lVert \mathcal{E}_A(\pi_1,\pi_2) \rVert = \lVert \pi_1 \rVert < \lVert \pi_1 \rVert + \lVert \pi_2 \rVert$.

Suppose the formula $A$ is an axiomatic formula in the final inference of $\pi_1$. Then the proof $\pi_1$ consists of an initial sequent and the multiset $\Gamma$ has the form $ \Gamma_0, A$. We obtain $\mathcal{E}_A(\pi_1,\pi_2)$ by applying the admissible rule $(\mathsf{weak})$:
\begin{gather*}
\AXC{$\pi_2$}
\noLine
\UIC{\vdots}
\noLine
\UIC{$A, \Sigma \Rightarrow \Pi$}
\LeftLabel{$(\mathsf{weak})$}
\RightLabel{ .}
\UIC{$\Gamma_0, A, \Sigma \Rightarrow  \Pi, \Delta$}
\DisplayProof
\end{gather*}
We have $\lVert \mathcal{E}_A(\pi_1,\pi_2) \rVert \leqslant \lVert \pi_2 \rVert  <\lVert \pi_1 \rVert + \lVert \pi_2 \rVert$.

Now suppose the formula $A$ is a side formula. Then the final inference in $\pi_1$ can be $(\mathsf{fix_L})$, $(\mathsf{fix_R})$, $(\mathsf{\rightarrow_L})$ or $(\mathsf{\rightarrow_R})$.

In the case of $(\mathsf{\rightarrow_R})$, the proof $\pi_1$ has the form
\begin{gather*}
\AXC{$\pi^\prime_1$}
\noLine
\UIC{\vdots}
\noLine
\UIC{$\Gamma, B \Rightarrow  C, \Delta_0, A$}
\LeftLabel{$(\mathsf{\rightarrow_R})$}
\RightLabel{ ,}
\UIC{$\Gamma \Rightarrow    B\rightarrow C, \Delta_0, A$}
\DisplayProof
\end{gather*}
where $\Delta = B\rightarrow C,\Delta_0$. We define $\mathcal{E}_A(\pi_1,\pi_2)$ as
\begin{gather*}
\AXC{$\mathcal{E}_A(\pi^\prime_1,\pi_2)$}
\noLine
\UIC{\vdots}
\noLine
\UIC{$\Gamma, B, \Sigma \Rightarrow  \Pi,C, \Delta_0$}
\LeftLabel{$(\mathsf{\rightarrow_R})$}
\RightLabel{ .}
\UIC{$\Gamma, \Sigma \Rightarrow  \Pi,  B\rightarrow C, \Delta_0$}
\DisplayProof
\end{gather*}
The proof $\mathcal{E}_A(\pi^\prime_1,\pi_2)$ is defined by the induction hypothesis for $\pi^\prime_1$ and $\pi_2$. We also have $\lVert \mathcal{E}_A(\pi_1,\pi_2) \rVert = \lVert \mathcal{E}_A(\pi^\prime_1,\pi_2) \rVert +1 < \lVert \pi^\prime_1 \rVert+ \lVert \pi_2 \rVert + 1= \lVert \pi_1 \rVert+ \lVert \pi_2 \rVert$.

In the case of $(\mathsf{fix_R})$, the proof $\pi_1$ has the form
\begin{gather*}
\AXC{$\pi^\prime_1$}
\noLine
\UIC{\vdots}
\noLine
\UIC{$\Gamma \Rightarrow  B[\mathsf{fp} \, x.\, B/\!/x] \Delta_0, A$}
\LeftLabel{$(\mathsf{fix_R})$}
\RightLabel{ ,}
\UIC{$\Gamma \Rightarrow  \mathsf{fp} \, x.\, B, \Delta_0, A$}
\DisplayProof
\end{gather*}
where $\Delta = \mathsf{fp} \, x.\, B,\Delta_0$. We define $\mathcal{E}_A(\pi_1,\pi_2)$ as
\begin{gather*}
\AXC{$\mathcal{E}_A(\pi^\prime_1,\pi_2)$}
\noLine
\UIC{\vdots}
\noLine
\UIC{$\Gamma,  \Sigma \Rightarrow  \Pi,B[\mathsf{fp} \, x.\, B/\!/p], \Delta_0$}
\LeftLabel{$(\mathsf{fix_R})$}
\RightLabel{ .}
\UIC{$\Gamma, \Sigma \Rightarrow  \Pi,  \mathsf{fp} \, x.\, B, \Delta_0$}
\DisplayProof
\end{gather*}
The proof $\mathcal{E}_A(\pi^\prime_1,\pi_2)$ is defined by the induction hypothesis, and $\lVert \mathcal{E}_A(\pi_1,\pi_2) \rVert = \lVert \mathcal{E}_A(\pi^\prime_1,\pi_2) \rVert +1 < \lVert \pi^\prime_1 \rVert+ \lVert \pi_2 \rVert +1= \lVert \pi_1 \rVert+ \lVert \pi_2 \rVert$.

The remaining cases of $(\mathsf{\rightarrow_L})$ and $(\mathsf{fix_L})$ can be analyzed analogously, so we omit them.

Now consider the final inference in $\pi_2$. If the formula $A$ is a weakening, an axiomatic or a side formula in it, then we can define $\mathcal{E}_A(\pi_1,\pi_2)$ in a similar way to the previous cases.

Suppose that the formula $A$ is a principal or an active formula in the final inferences of $\pi_1$ and $\pi_2$.
Then $A$ has the form $\mathsf{fp} \, x. \,A_0$, $A_0\rightarrow A_1$ or $\Box A_0$.

If $A=\Box A_0$, then
$\pi_2$ has one of the two forms
\begin{gather*}
\AXC{$\pi^\prime_2$}
\noLine
\UIC{\vdots}
\noLine
\UIC{$ A_0, \Sigma_1, \Box \Sigma_2\Rightarrow  D$}
\LeftLabel{$(\mathsf{\Box})$}
\RightLabel{ }
\UIC{$ \Sigma_0, \Box A_0 ,\Box \Sigma_1, \Box \Sigma_2 \Rightarrow   \Box D, \Pi_0$}
\DisplayProof \qquad
\AXC{$\pi^\prime_2$}
\noLine
\UIC{\vdots}
\noLine
\UIC{$  \Sigma_1, \Box A_0,\Box \Sigma_2\Rightarrow  D$}
\LeftLabel{$(\mathsf{\Box})$}
\RightLabel{ ,}
\UIC{$ \Sigma_0, \Box \Sigma_1, \Box A_0 , \Box \Sigma_2  \Rightarrow   \Box D, \Pi_0$}
\DisplayProof
\end{gather*}
where $\Sigma = \Sigma_0, \Box \Sigma_1, \Box \Sigma_2$ and $\Pi = \Box D, \Pi_0$.
In addition, the proof $\pi_1$ has the form
\begin{gather*}
\AXC{$\pi^\prime_1$}
\noLine
\UIC{\vdots}
\noLine
\UIC{$\Gamma_1, \Box \Gamma_2 \Rightarrow   A_0$}
\LeftLabel{$(\mathsf{\Box})$}
\RightLabel{ ,}
\UIC{$\Gamma_0, \Box \Gamma_1, \Box \Gamma_2 \Rightarrow  \Box A_0,   \Delta$}
\DisplayProof
\end{gather*}
where $\Gamma = \Gamma_0, \Box \Gamma_1, \Box \Gamma_2$. If $\pi_2$ has the first form, then we define $\mathcal{E}_A(\pi_1,\pi_2)$ as
\begin{gather*}
\AXC{$\mathcal{E}_{A_0}(\pi^\prime_1,\pi^\prime_2)$}
\noLine
\UIC{\vdots}
\noLine
\UIC{$\Gamma_1, \Box \Gamma_2, \Sigma_1, \Box \Sigma_2 \Rightarrow  D$}
\LeftLabel{$(\mathsf{\Box})$}
\RightLabel{ .}
\UIC{$\Gamma_0, \Box \Gamma_1, \Box \Gamma_2,\Sigma_0, \Box \Sigma_1, \Box \Sigma_2 \Rightarrow   \Box D, \Pi_0, \Delta$}
\DisplayProof
\end{gather*}
We have $\lVert \mathcal{E}_A(\pi_1,\pi_2) \rVert = \lVert \mathcal{E}_{A_0}(\pi^\prime_1,\pi^\prime_2)\rVert +1 < \lVert \pi^\prime_1\rVert + \lVert \pi^\prime_2\rVert +1 < \lVert \pi_1\rVert + \lVert \pi_2\rVert$.
If $\pi_2$ has the second form, then we define $\mathcal{E}_A(\pi_1,\pi_2)$ as
\begin{gather*}
\AXC{$\mathcal{E}_A(f(\pi_1),\pi^\prime_2)$}
\noLine
\UIC{\vdots}
\noLine
\UIC{$\Box \Gamma_1, \Box \Gamma_2, \Sigma_1, \Box \Sigma_2 \Rightarrow  D$}
\LeftLabel{$(\mathsf{\Box})$}
\RightLabel{ ,}
\UIC{$\Gamma_0, \Box \Gamma_1, \Box \Gamma_2,\Sigma_0, \Box \Sigma_1, \Box \Sigma_2 \Rightarrow   \Box D, \Pi_0, \Delta$}
\DisplayProof
\end{gather*}
where $f(\pi_1)$ is the proof obtained by erasing multisets $\Gamma_0$ and $\Delta$ from the conclusion of $\pi_1$.
We have $\lVert \mathcal{E}_A(\pi_1,\pi_2) \rVert = \lVert \mathcal{E}_A(f(\pi_1),\pi^\prime_2)\rVert +1 < \lVert f(\pi_1)\rVert + \lVert \pi^\prime_2\rVert +1 = \lVert \pi_1\rVert + \lVert \pi_2\rVert$.

In the case of $A = \mathsf{fp} \, x. \,A_0 $, the proofs $\pi_1$ and $\pi_2$ have the form
\begin{gather*}
\AXC{$\pi^\prime_1$}
\noLine
\UIC{\vdots}
\noLine
\UIC{$\Gamma \Rightarrow  \Delta, A_0[ \mathsf{fp} \, x. \,A_0 /\!/x]$}
\LeftLabel{$(\mathsf{fix_R})$}
\UIC{$\Gamma\Rightarrow  \Delta,  \mathsf{fp} \, x. \,A_0 $}
\DisplayProof \qquad
\AXC{$\pi^\prime_2$}
\noLine
\UIC{\vdots}
\noLine
\UIC{$A_0[ \mathsf{fp} \, x. \,A_0 /\!/x], \Sigma \Rightarrow  \Pi$}
\LeftLabel{$(\mathsf{fix_L})$}
\RightLabel{ .}
\UIC{$\mathsf{fp} \, x. \,A_0  ,\Sigma \Rightarrow  \Pi$}
\DisplayProof
\end{gather*}
We put $\mathcal{E}_A(\pi_1,\pi_2) = \mathcal{E}_{A_0[ \mathsf{fp} \, x. \,A_0 /\!/x]}(\pi^\prime_1,\pi^\prime_2)$ and see that
$\lVert \mathcal{E}_A(\pi_1,\pi_2)\rVert
  = \lVert \mathcal{E}_{A_0[ \mathsf{fp} \, x. \,A_0 /\!/x]}(\pi^\prime_1,\pi^\prime_2)\rVert
   < \lVert \pi^\prime_1 \rVert + \lVert \pi^\prime_2 \rVert < \lVert \pi_1
    \rVert + \lVert\pi_2 \rVert$.

If $A=A_0 \rightarrow A_1$, then the proofs $\pi_1$ and $\pi_2$ have the form
\begin{gather*}
\AXC{$\pi^\prime_1$}
\noLine
\UIC{\vdots}
\noLine
\UIC{$A_0, \Gamma \Rightarrow  \Delta, A_1$}
\LeftLabel{$(\mathsf{\rightarrow_R})$}
\UIC{$\Gamma\Rightarrow  \Delta, A_0 \rightarrow A_1$}
\DisplayProof \qquad
\AXC{$\pi^\prime_2$}
\noLine
\UIC{\vdots}
\noLine
\UIC{$ A_1, \Sigma_1 \Rightarrow  \Pi_1$}
\AXC{$\pi^{\prime \prime}_2$}
\noLine
\UIC{\vdots}
\noLine
\UIC{$  \Sigma_0 \Rightarrow  \Pi_0, A_0$}
\LeftLabel{$(\mathsf{\rightarrow_L})$}
\RightLabel{ ,}
\BIC{$\Sigma_0, A_0 \rightarrow A_1 ,\Sigma_1 \Rightarrow  \Pi_0, \Pi_1$}
\DisplayProof
\end{gather*}
where $\Sigma =\Sigma_0, \Sigma_1$ and $\Pi= \Pi_0, \Pi_1 $. By the induction hypothesis, $\mathcal{E}_{A_0}(\pi^{\prime\prime}_2,\pi^\prime_1)$ is defined and $ \lVert \mathcal{E}_{A_0}(\pi^{\prime\prime}_2,\pi^\prime_1)\rVert < \lVert\pi^{\prime\prime}_2 \rVert +\lVert\pi^\prime_1 \rVert$.
Since $\lVert\mathcal{E}_{A_0}(\pi^{\prime\prime}_2,\pi^\prime_1) \rVert+ \lVert\pi^\prime_2 \rVert<\lVert\pi^{\prime\prime}_2 \rVert +\lVert\pi^\prime_1 \rVert+ \lVert\pi^\prime_2 \rVert<\lVert\pi_1 \rVert+ \lVert\pi_2 \rVert $, then $ \mathcal{E}_{A_1} (\mathcal{E}_{A_0}(\pi^{\prime\prime}_2,\pi^\prime_1), \pi^\prime_2)$ is defined by the induction hypothesis. We put $\mathcal{E}_A(\pi_1,\pi_2) = \mathcal{E}_{A_1} (\mathcal{E}_{A_0}(\pi^{\prime\prime}_2,\pi^\prime_1), \pi^\prime_2)$. In addition, we have
$\lVert \mathcal{E}_A(\pi_1,\pi_2)\rVert  = \lVert \mathcal{E}_{A_1} (\mathcal{E}_{A_0}(\pi^{\prime\prime}_2,\pi^\prime_1), \pi^\prime_2)\rVert < \lVert\mathcal{E}_{A_0}(\pi^{\prime\prime}_2,\pi^\prime_1) \rVert+ \lVert\pi^\prime_2 \rVert<\lVert\pi_1 \rVert+ \lVert\pi_2 \rVert$.
\end{proof}

\section{Conclusions and future work}

The preliminary results presented in this paper indicate the following conclusions:

\bi
\item Derivability conditions can be stated in a way not assuming much about logic. However,
\item G\"odel's argument presupposes a certain amount of contraction for the logic under consideration.
\ei
The role of contraction rule here is somewhat similar to its role in Liar-type paradoxes including Russell's paradox in set theory. Thus, Vyacheslav Grishin (see~\cite{Gri74,Gri82}) pioneered the study of set theory with full comprehension based on a logic without contraction. He demonstrated that the pure comprehension scheme is consistent in this logic. He also showed, however, that the extensionality principle allows for this system to actually \emph{prove} contraction even if there is no postulated contraction in the logic.

One can also consider systems of arithmetic based on contraction-free logic, see e.g.~Restall~\cite[Chapter 11]{Res94}. For one such system, considered by the second author of this paper, the rule of $\Box$-contraction is admissible, which according to our results still yields G2. Thus, we are still missing convincing examples of mathematical theories based on weak logics for which G2 would fail.

\bi
\item For consequence relations with an implication and with $\Box$ satisfying L\"ob's conditions, the existence of appropriately many fixed points does not imply their uniqueness. Nor does it imply formalized versions of G2 and L\"ob's theorem $\Box(\Box\phi\to\phi)\vdash \Box\phi$.
\ei
This shows that the move from \emph{diagonalized algebras} in the sense of R.~Magari, i.e., Boolean algebras with $\Box$ satisfying L\"ob's conditions and having fixed points, to \emph{diagonalizable algebras} (modal algebras satisfying L\"ob's identity) is, in general, not possible for logics without contraction and weakening. See Smory\'nski~\cite{Smo82a,Smo85} for a nice exposition of the original setup.

\bi
\item One can also show that the admissibility of L\"ob's rule does not, in general, imply a formalized version of G2.
\ei
A system $\mathsf S^*$ witnessing this property can be obtained by extending the notion of proof in the system $\mathsf S$ to possibly non-well-founded proof trees. Infinite proofs may arise because of the presence of the fixed point rules. For $\mathsf S^*$, unlike $\mathsf S$, one can show that L\"ob's rule is admissible. Yet, formalized G2 is still underivable.
The analysis of $\mathsf S^*$ is based on another cut-admissibility theorem, which we postpone to a later publication.

We remark that the system $\mathsf S$ does not provide a counterexample to the non-formalized version of G2, since $\Imp \neg \Box\bot$ is not provable. We believe that such a counterexample can be constructed by extending the language of $\mathsf S$ by an operator similar to $!$ from linear logic and adding to $S$ a fixed point of the form $a=\Diamond!a$. However, a confirmation of this hypothesis is left for future work.

\section{Acknowledgements}
The authors would like to thank Johan van Benthem for useful comments and questions. This work is supported by the Russian Foundation for Basic Research, grant 15-01-09218a, and by the Presidential council for support of leading scientific schools.


\input{gen_goedel.bbl}
\end{document}

%% file: macros.tex
\newtheorem{theorem}{Theorem}
\newtheorem{lemma}{Lemma}[section]

\newtheorem{corollary}[lemma]{Corollary}
\newtheorem{proposition}[lemma]{Proposition}

\newtheorem{definition}[lemma]{Definition}
\newtheorem{example}[lemma]{Example}
\newtheorem{remark}[lemma]{Remark}

\newcommand{\bl}{\begin{lemma}}
\newcommand{\el}{\end{lemma}}
\newcommand{\bt}{\begin{theorem}}
\newcommand{\et}{\end{theorem}}
\newcommand{\bcor}{\begin{corollary}}
\newcommand{\ecor}{\end{corollary}}
\newcommand{\bp}{\begin{proof}}
\newcommand{\ep}{\end{proof}}
\newcommand{\bpr}{\begin{proposition}}
\newcommand{\epr}{\end{proposition}}
\newcommand{\brem}{\begin{remark} \em}
\newcommand{\erem}{\end{remark}}
\newcommand{\bd}{\begin{definition} \em}
\newcommand{\ed}{\end{definition}}
\newcommand{\bex}{\begin{example} \em
}
\newcommand{\eex}{\end{example}}
\newcommand{\beq}{\begin{equation} }
\newcommand{\eeq}{\end{equation}}

\newcommand{\bi}{\begin{itemize}
  }
\newcommand{\ei}{\end{itemize}}
\newcommand{\ben}{\begin{enumerate} }
\newcommand{\een}{\end{enumerate} }

\newenvironment{enumr}{

\renewcommand{\labelenumi}{$\mathrm{(\theenumi)}$}
\begin{enumerate}     }{\end{enumerate}

\renewcommand{\labelenumi}{\theenumi.}}

\newcommand{\al}[1]{\forall #1\:}

\newcommand{\ga}{\alpha}

\renewcommand{\phi}{\varphi}

\newcommand{\Prs}[1]{\mathsf{Prov}_{#1}}

\renewcommand{\Pr}{\mathrm{Pr}}

\newcommand{\Con}{\mathrm{Con}}

\newcommand{\PA}{\mathsf{PA}}

\newcommand{\Imp}{\Rightarrow}

\newcommand{\Boxt}{\boxtimes}
\renewcommand{\Box}{\square}

\renewcommand{\leq}{\leqslant}

\newcommand{\eop}{$\dashv$ \protect\par \addvspace{\topsep}}


%% file: gen_goedel_all1.bbl
\begin{thebibliography}{10}

\bibitem{AZ11}
Z.~Adamowicz and K.~Zdanowski.
\newblock Lower bounds for the provability of herbrand consistency in weak
  arithmetics.
\newblock {\em Fundamenta Mathematicae}, 212(3):191--216, 2011.

\bibitem{BezShe}
A.~Bezboruah and J.~C. Shepherdson.
\newblock G\"{o}del's second incompleteness theorem for {Q}.
\newblock {\em The Journal of Symbolic Logic}, 41(2):503--512, 1976.

\bibitem{Fef60}
S.~Feferman.
\newblock Arithmetization of metamathematics in a general setting.
\newblock {\em Fundamenta Mathematicae}, 49:35--92, 1960.

\bibitem{God-CW-I}
S.~Feferman, J.R. Dawson, S.C. Kleene, G.H. Moore, R.M. Solovay, and J.~van
  Heijenoort, editors.
\newblock {\em Kurt G\"{o}del Collected Works, Volume 1: Publications
  1929--1936}. Oxford Univeristy Press, 1996.

\bibitem{God31}
K.~G\"{o}del.
\newblock \"{U}ber formal unentscheidbare {S\"atze} der {Principia Mathematica}
  und verwandter {Systeme I}.
\newblock {\em Monatshefte f\"ur Mathematik und Physik}, 38:173--198, 1931.

\bibitem{God33}
K.~G\"{o}del.
\newblock Eine {I}nterpretation des intuitionistischen {A}ussagenkalkuls.
\newblock {\em Ergebnisse Math. Kolloq.}, 4:39--40, 1933.
\newblock English translation in \cite{God-CW-I}, pages 301--303.

\bibitem{Gri74}
V.N. Grishin.
\newblock On some non-standard logic and its application to set theory.
\newblock In {\em Investigations on formalized languages and non-classical
  logics}, pages 135--171. Nauka, Moscow, 1974.
\newblock In Russian.

\bibitem{Gri82}
V.N. Grishin.
\newblock Predicate and set-theoretic calculi based on logic without
  contractions.
\newblock {\em Mathematics of the USSR-Izvestiya}, 18(1):41--59, 1982.

\bibitem{HB}
D.~Hilbert and P.~Bernays.
\newblock {\em Grundlagen der Mathematik, Vols. I and II, 2d ed.}
\newblock Springer-{V}erlag, Berlin, 1968.

\bibitem{Jer70}
R.G. Jeroslow.
\newblock Consistency statements in formal theories.
\newblock {\em Fundamenta Mathematicae}, 72:2--39, 1970.

\bibitem{Jer73}
R.G. Jeroslow.
\newblock Redundancies in the {Hilbert--Bernays} derivability conditions.
\newblock {\em The Journal of Symbolic Logic}, 38(3):359--367, 1973.

\bibitem{Kle50}
S.C. Kleene.
\newblock A {symmetric form of G\"odel's theorem}.
\newblock {\em {Indagationes Mathematicae}}, 12:244--246, 1950.

\bibitem{Lob55}
M.H. L\"{o}b.
\newblock Solution of a problem of {Leon Henkin}.
\newblock {\em The Journal of Symbolic Logic}, 20:115--118, 1955.

\bibitem{Pud85}
P.~Pudl\'ak.
\newblock Cuts, consistency statements and interpretations.
\newblock {\em The Journal of Symbolic Logic}, 50:423--441, 1985.

\bibitem{Rau06}
W.~Rautenberg.
\newblock {\em A Concise Introduction to Mathematical Logic}.
\newblock Springer, second edition, 2006.

\bibitem{Res94}
G.~Restall.
\newblock {\em On Logics Without Contraction}.
\newblock PhD thesis, The University of Queensland, 1994.
\newblock http://consequently.org/papers/onlogics.pdf.

\bibitem{Smo77}
C.~Smory\'{n}ski.
\newblock {T}he incompleteness theorems.
\newblock In J.~Barwise, editor, {\em {H}andbook of {M}athematical {L}ogic},
  pages 821--865. North Holland, Amsterdam, 1977.

\bibitem{Smo82a}
C.~Smory\'{n}ski.
\newblock Fixed point algebras.
\newblock {\em Bull. Amer. Math. Soc.}, 6(3):317--356, 1982.

\bibitem{Smo85}
C.~Smory\'{n}ski.
\newblock {\em Self-Reference and Modal Logic}.
\newblock Springer-Verlag, Berlin, 1985.

\bibitem{Smu94}
R.M. Smullyan.
\newblock {\em Diagonalization and Self-Reference}.
\newblock Oxford Logic Guides 27. Oxford University Press, 1994.

\bibitem{Vis93}
A.~Visser.
\newblock Unprovability of small inconsistency.
\newblock {\em Archive for Math.\ Logic}, 32:275--298, 1993.

\bibitem{Vis11}
A.~Visser.
\newblock Can we make the {Second Incompleteness Theorem} coordinate free?
\newblock {\em Journal of Logic and Computation}, 21(4):543--560, 2011.

\bibitem{Vis12}
A.~Visser.
\newblock The {Second Incompleteness Theorem} and bounded interpretations.
\newblock {\em Studia Logica}, 100(1--2):399--418, 2012.

\bibitem{WP}
A.~Wilkie and J.~Paris.
\newblock On the scheme of induction for bounded arithmetic formulas.
\newblock {\em Annals of Pure and Applied Logic}, 35:261--302, 1987.

\bibitem{Wil01}
D.~Willard.
\newblock Self-verifying systems, the incompleteness theorem and the
  tangibility reflection principle.
\newblock {\em The Journal of Symbolic Logic}, 66:536--596, 2001.

\bibitem{Wil06}
D.~Willard.
\newblock A generalization of the {Second Incompleteness Theorem} and some
  exceptions to it.
\newblock {\em Annals of Pure and Applied Logic}, 141:472--496, 2006.

\end{thebibliography}
